\documentclass[12pt]{amsart}

\usepackage{amsthm}
\usepackage{epsfig}

\newtheorem{thm}{Theorem}
\newtheorem{lemma}[thm]{Lemma}

\newtheorem{prop}[thm]{Proposition}
\theoremstyle{definition}
\newtheorem{defn}{Definition}
\newtheorem*{notn}{Notation}

\def\interval#1#2{I_{#1}(#2)}
\def\Itop#1{\interval{top}{#1}}
\def\Ibot#1{\interval{bot}{#1}}

\def\fig#1#2#3#4{
\begin{figure}[ht]
\begin{center}
\epsfig{file=#1,width=#3}
\caption{#2}
\label{#4}
\end{center}
\end{figure}
}

\def\fign#1#2#3{
\begin{figure}[ht]
\begin{center}
\epsfig{file=#1,width=#3}
\caption{#2}
\end{center}
\end{figure}
}

\def\arlt{\overset{\downarrow}{<}}
\def\arlte{\overset{\downarrow}{\leq}}

\begin{document}

\title{Proper and Unit Trapezoid Orders and Graphs}

\author{Kenneth P. Bogart}
\address{\hskip -\parindent
Kenneth P. Bogart\\
Department of Mathematics
\\6188 Bradley Hall
\\Dartmouth College
\\Hanover, NH 03755}
\email{kenneth.p.bogart@dartmouth.edu}

\author{Rolf H. M\"ohring}
\address{\hskip -\parindent
Rolf H. M\"ohring
\\Fachbereich Mathematik
\\Sekr. 6-1\\
Technische Universit{\"a}t Berlin\\
Stra{\ss}e des 17. Juni 136\\
10623  Berlin, Germany}
\email{moehring@math.tu-berlin.de}
\author{Stephen P. Ryan}

\address{\hskip-\parindent
Stephen Ryan
\\Department of Mathematics\\
6188 Bradley Hall\\
Dartmouth College
\\Hanover, NH 03755}

\email{stephen.p.ryan@dartmouth.edu}

\thanks{Supported by ONR contract N0014-91-J1019} 
\thanks{Supported by Dartmouth College as Harris German-Dartmouth
Distinguished Visiting Professor}
\thanks{Supported by ONR contract N0014-94-1-0950}
\thanks{Research at MSRI is supported in part by NSF grant
DMS-9022140}

\date{18 November 1996}

\begin{abstract}We show that the class of trapezoid orders in which no
trapezoid strictly contains any other trapezoid strictly contains the
class of trapezoid orders in which every trapezoid can be drawn with
unit area.  This is different from the case of interval orders, where
the class of proper interval orders is exactly the same as the class
of unit interval orders.
\end{abstract}

\maketitle

\par\bigskip

\section*{Introduction}

\begin{defn}
An order $\prec$ of a set $X$ is called a {\bf trapezoid order} if,
given two parallel lines in the plane (which we will take to be
horizontal), there is, associated to each $x \in X$, a trapezoid $T_x$
with one base on each of the lines with the property that $x\prec y$
if and only if $T_x \cap T_y = \emptyset$ and each point in $T_x$ is
to the left of some point in $T_y$. \cite{langley:tolerance}
\end{defn}

A trapezoid order is a generalization of an interval order, in which
there is an interval $I_x$ associated to each element $x$ of the
order, with $x \prec y$ if and only if each point of $I_x$ lies to the
left of each point of $I_y$.  We will only deal with interval orders
briefly, so the interested reader should see \cite{fishburn85} or
\cite{golumbic80:_algorithmic} for a treatment of interval orders.
Interval orders first appeared in \cite{wiener14}
\cite{fishburn92:_norbert_wiener}, but were not studied again until
\cite{fishburn70:_intran} and \cite{mirkin72:_descr}.

\begin{notn}
Such a collection of trapezoids will be called a {\bf trapezoid
representation} of the order.  Thus, an order is a trapezoid order if
and only if it has a trapezoid representation.
\end{notn}

\begin{notn}
The two parallel lines will be called the baselines, from the fact
that all the trapezoids have their bases on those two lines.  
\end{notn}

It may at times be convenient to confuse the elements of
the ordered set with the trapezoids representing those elements.

\begin{notn}
Each trapezoid defines two intervals, one on each baseline, determined
by the intersection of the trapezoid with the baseline.  We call the
interval determined in this way by the intersection of the trapezoid
with the upper baseline the {\bf upper interval}, or {\bf top}
interval of $T_x$ and refer to it by $\Itop{x} = [L(x),R(x)]$.  Similarly, the
interval determined by intersection with the lower baseline is called
the {\bf lower interval}, or {\bf bottom} interval of $T_x$ and
referred to by $\Ibot{x} = [l(x),r(x)]$.
\end{notn}

\begin{notn}
The length of the upper interval will be denoted by $t(x)$; it is, of
course, just $R(x)-L(x)$, but we would like to avoid excess notation
later on.  Similarly, the length of the lower interval will be denoted
by $b(x)$.
\end{notn}

Conversely, given two intervals on parallel lines, we can define a
trapezoid by taking the convex hull of the two intervals.  This
defines an immediate and obvious bijection between trapezoids on the
two parallel lines, and pairs of intervals, one on each line.

The collection of all the upper intervals $\{\Itop{x}|x \in X\}$
defined as above gives rise to another ordering on the original set, a
natural interval ordering.  Similarly, there is an interval ordering
defined on the set by the collection of lower intervals.  It is easy
to check that $x \prec y$ (in the trapezoid ordering) if and only if
$\Itop{x} \prec \Itop{y}$ and $\Ibot{x} \prec \Ibot{y}$ (in the
respective interval orderings).  However, the latter condition is
precisely the definition of the intersection of the two interval
orderings.  Hence, a trapezoid order is the intersection of two
interval orders and therefore has interval dimension at most two.

Conversely, given an order with interval dimension at most two, we can
take any two (possibly identical) interval orders whose intersection
is the order under consideration, place interval representations for
each on two parallel lines, and use the bijection discussed above to
produce trapezoids associated to each element.  The observation that
$x \prec y$ in the trapezoid ordering precisely when $x \prec y$ in
each of the two interval orderings assures us that the result is a
trapezoid representation of the given order, and thus that the order
is a trapezoid order.

Thus, the class of trapezoid orders is precisely the class of orders
whose interval dimension is at most two \cite{langley:tolerance}.

Note that a trapezoid order may have (in fact, must have) many
trapezoid representations, since rescalings and small perturbations in
the endpoints do not change the underlying order.

We will use the ``$\prec$'' symbol to indicate the predecessor
relation in ordered sets, and we will use ``$\|$'' to refer to
incomparability in ordered sets.  We will also make use of the
ordering of the real numbers, when referring to endpoints and
intervals defined by the trapezoids of trapezoid orders, and use the
standard ``$<$'' symbol for that ordering.

\begin{defn}A {\bf proper} trapezoid order is one for which there is a 
trapezoid representation in which no trapezoid is properly contained
in any other.
\end{defn}

\begin{defn}A {\bf unit} trapezoid order is one for which there is a
trapezoid representation in which every trapezoid has the same area.  By
a suitable choice of scaling, we may assume that this area is 1 and
that the distance between the two baselines is also 1 (and hence, that
the sum of the lengths of the two bases of any trapezoid in a
representation of this order is 2).
\end{defn}

These definitions are motivated by corresponding definitions for
interval orders.  \cite{roberts69:indifference} proved, in the context
of interval graphs, that proper and unit interval orders are
equivalent.

\begin{defn} A graph $G=(V,E)$ is said to be a {\bf trapezoid graph}
if there are two parallel lines such that for each $x \in V$ there is
a trapezoid $T_x$ associated to it with one base on each line such
that $(x,y) \in E$ if and only if $T_x \cap T_y \not = \emptyset$.
\end{defn}

Trapezoid graphs were first discussed in \cite{dagan88:trapezoid}.

There is a corresponding notion of {\bf interval graphs}, introduced
by \cite{hajos57:_uber_art} and \cite{benzer59}, but first called
interval graphs by \cite{gilmore62}\nocite{gilmore64}.

Clearly, a graph is a trapezoid graph if and only if it is the
cocomparability graph of a trapezoid order (partial order having
interval dimension at most two) \cite{dagan88:trapezoid}.
This is also a consequence of \cite{habib91:_inter}.

The class of trapezoid orders (graphs) is also the same as the class
of bounded bitolerance orders (graphs), though the natural definitions
of proper and unit are different in that context
\cite{langley:tolerance}. Bogart and Isaak proved that proper and unit
in that context are equivalent \cite{bogart:bitolerance}. A similiar
result for digraphs was proved in \cite{shull95:_unit_proper}.

We will be approaching the main results from an order-theoretic point
of view.  By analyzing the autonomous sets of the orders, we shall
show that all of the following existence results apply to trapezoid
graphs, as well, so that there exists a trapezoid graph which is an
improper trapezoid graph, and there exists a trapezoid graph which is
a proper trapezoid graph, but not a unit trapezoid graph.

\section*{Preliminaries}

\fig{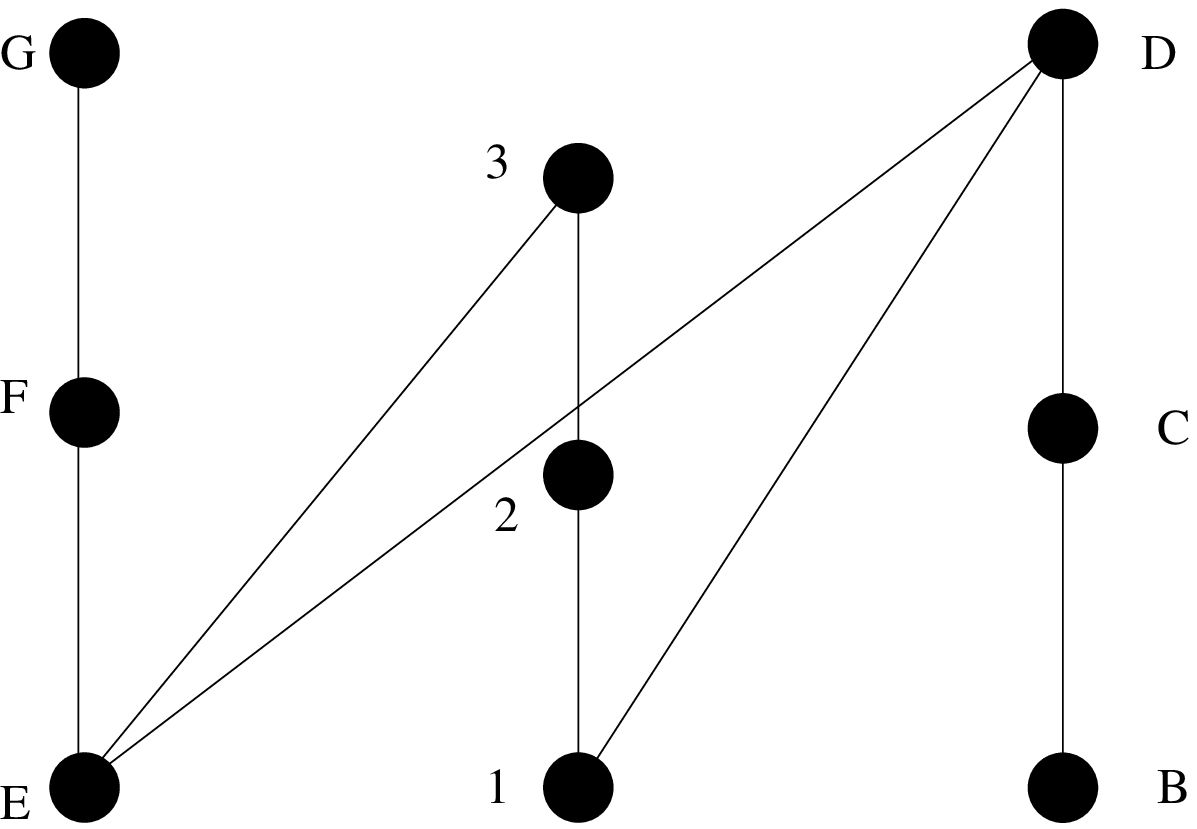}{The jaw order}{3in}{fig:jaworder}

\begin{lemma}[Jaw Lemma]  
\label{lem:jaw}

The order shown in Figure \ref{fig:jaworder} has a trapezoid
representation, and hence is a trapezoid order.  Further, in every
trapezoid representation, we must have endpoints in the relations
$$r(B)<l(C)\leq r(1)<l(2)\leq r(E)<l(D)\leq r(2)<l(3)\leq r(F)<l(G)$$ and
$$R(E)< L(2)\leq R(2)< L(D)$$or in the relations
$$R(B)<L(C)\leq R(1)<L(2)\leq R(E)<L(D)\leq R(2)<L(3)\leq R(F)<L(G)$$ and
$$r(E)< l(2)\leq r(2)< l(D)$$
\end{lemma}

\begin{proof}
First, a word on the name of the lemma.  Figure \ref{fig:jaws} is a
portion of a generic trapezoid representation of the order in Figure
\ref{fig:jaworder}.  The name comes from the fact that any element
placed above $E$ and below $D$ will be forced to have the interval
which is its base on this line completely contained in the interval
for $2$.  The appearance of the trapezoids for $D$ and $E$ on this
side is that of teeth, and their purpose is to squeeze other elements
completely inside the interval for $2$ on this side.

\fig{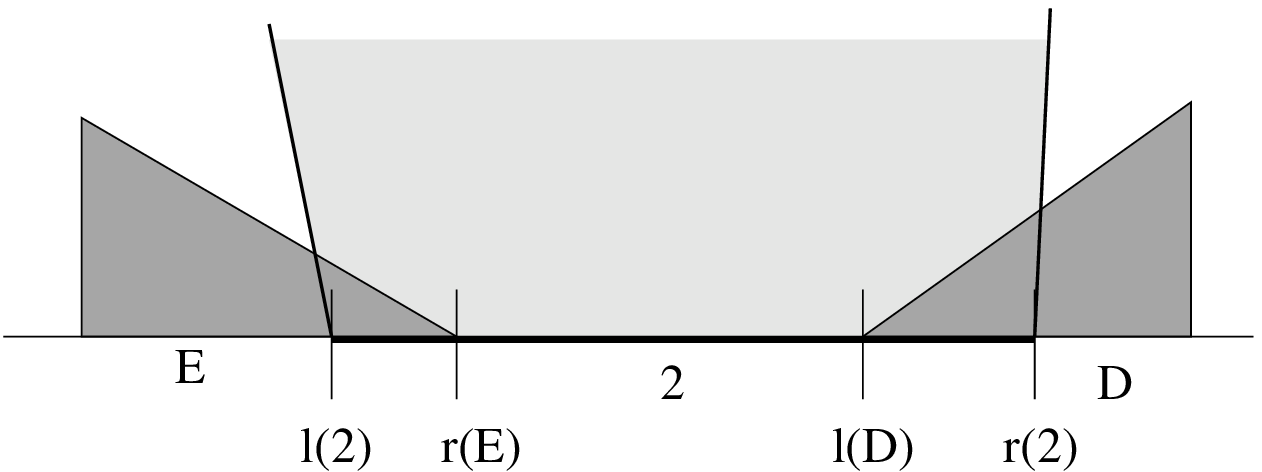}{The ``jaws'' of the jaw
lemma}{\hsize}{fig:jaws}

For the proof, we will build up a portion of a trapezoid
representation, using only those features that are common to all
trapezoid representations of that order.  There is some room for
variation on the fringes of the resulting diagram, but the center
portion of the diagram will be forced to appear as in Figure
\ref{fig:jaws}.

We first observe that $1 \prec 2 \prec 3$ is a chain, and hence will
have to take the form of three trapezoids in a row.  (See Figure
\ref{fig:chain}) \fig{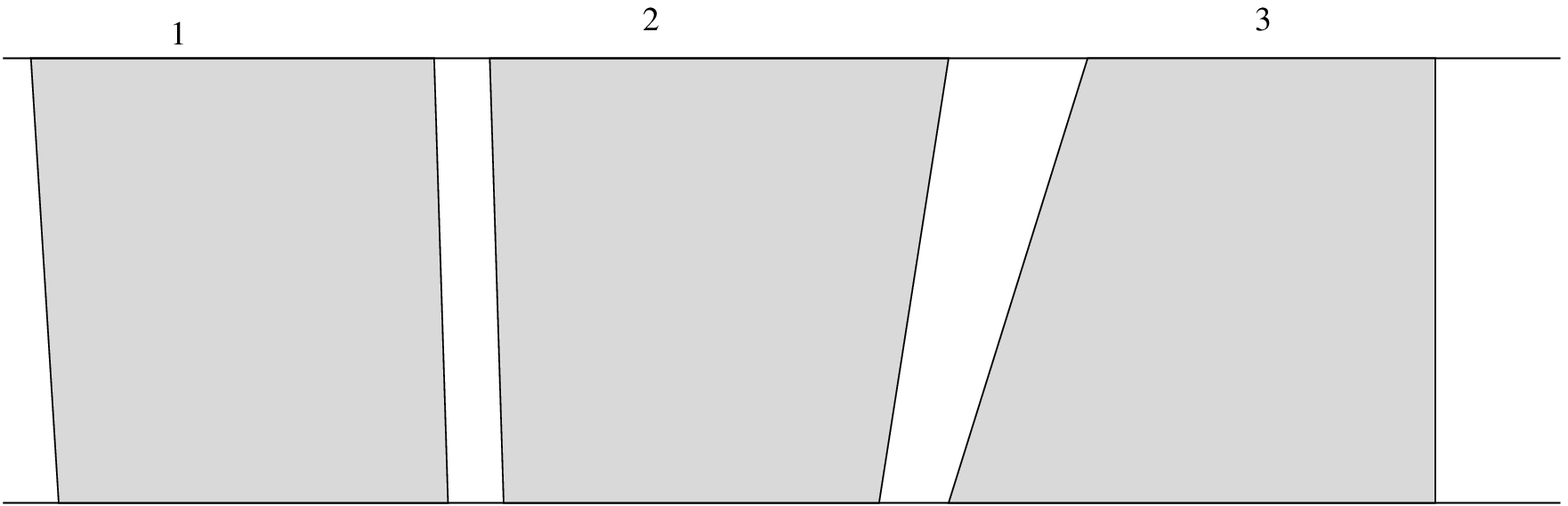}{The chain
trapezoids}{\hsize}{fig:chain} This chain is used as a sort of
coordinate system to aid in identifying (and forcing!) relative
positions of other trapezoids.  In all of the following exposition, it
will be assumed that $l(1) \leq r(1) < l(2) \leq r(2) < l(3) \leq
r(3)$ and $L(1) \leq R(1) < L(2) \leq R(2) < L(3) \leq R(3)$, even if
these relations are not explicitly stated again.

Corresponding inequalities, of course, are true for any chain.  We
shall prove the inequalities shown below.  The marked inequalities
follow from chains in Figure \ref{fig:jaworder}.  The others are to be
explained in what follows.
$$r(B)\arlt l(C)\leq r(1)\arlt l(2)\leq r(E)\arlt l(D)\leq r(2)\arlt
l(3)\leq r(F)\arlt l(G)$$
$$R(E)<L(2)\arlte R(2)< L(D)$$

$D$ is incomparable to 2 and 3, but is over 1.  Therefore, $T_D$ must
lie entirely to the right of $T_1$ while still overlapping $T_2$ and
$T_3$.  Since $T_D$ overlaps $T_2$, we may conclude that either $l(D) \leq
r(2)$ or that $L(D) \leq R(2)$.  Without loss of generality, assume the
lower baseline, i.e.
\begin{equation}\label{eq:assumption}l(D) \leq r(2)\end{equation} If it
is the upper baseline, an entirely symmetric argument will prove the
other statement in the lemma.  In our chain of inequalities, this
proves
$$ r(B)\arlt l(C)\leq r(1)\arlt l(2)\leq r(E)\arlt l(D)
\overset{\eqref{eq:assumption}}{\arlte} r(2)\arlt l(3)\leq r(F)\arlt l(G)$$

Now $B$ is below $D$, yet incomparable to 3.  Since $l(D) \leq r(2)$, it
must be that $T_B$ cannot overlap $T_3$ on the lower baseline, thus it
must overlap $T_3$ on the upper baseline; i.e. $R(B) \geq L(3)$.  This is
clearer when thought of as the intersection of two interval orders: on
the lower baseline, $\Ibot{B} \prec \Ibot{3}$ because of the
positioning of the left endpoint of $\Ibot{D}$; since $B \| 3$, it
must be the case that on the upper baseline, the interval for $B$
overlaps the interval for 3, or lies completely to the right of it; at
the moment, we don't care, just so long as $\Itop{B} \not\prec
\Itop{3}$.  In either case, $L(3)\leq R(B)$.  However, from the chains in
Figure \ref{fig:jaworder},
\begin{equation}\label{eq:incomp1}R(2)<L(3)\leq R(B)<L(D)\end{equation}
giving us
$$R(E)<L(2)\arlte R(2)\overset{\eqref{eq:incomp1}}{\arlt} L(D)$$

Additionally, since $C\|1$, it follows that $l(C)\leq r(1)$, by a
symmetric argument to the one in the previous paragraph.  From
\eqref{eq:incomp1}, $L(3)\leq R(B)$, while from the chains of Figure
\ref{fig:jaworder}, $R(B)<L(C)$, so in order for $T_C$ to overlap
$T_1$, it must be the case that $\Ibot{C}$ lies to the left of or
overlaps $\Ibot{1}$.  In symbols,
\begin{equation}\label{eq:Cleft}l(C)\leq r(1)\end{equation} which gives us
$$r(B)\arlt l(C)\overset{\eqref{eq:Cleft}}{\arlte} r(1)\arlt
l(2)\leq r(E)\arlt l(D)\overset{\eqref{eq:assumption}}{\arlte} r(2)\arlt
l(3)\leq r(F)\arlt l(G)$$

\fig{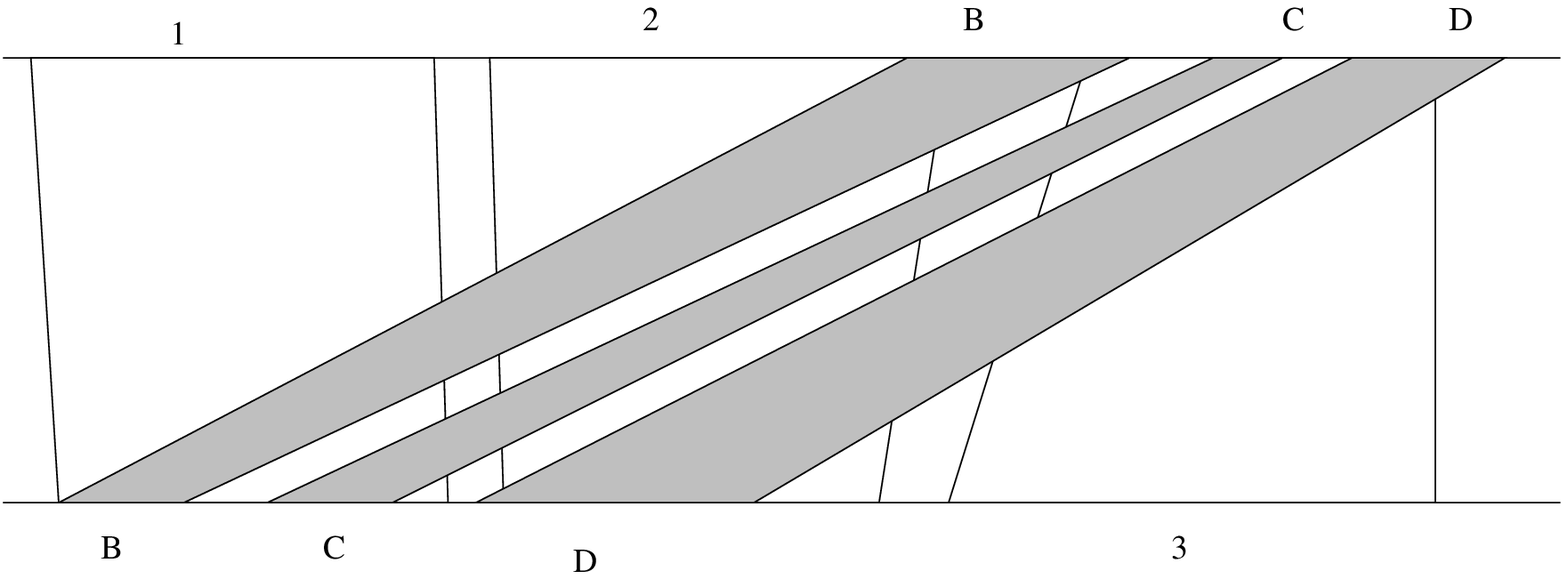}{A partial trapezoid representation of the jaw order,
illustrating the inequalities proved through Equation \eqref{eq:Cleft}.  
According to what has been proved so far, we have no information as to
the relative positions of $l(D)$ and $l(2)$, nor of several other
pairs of endpoints; this is, however, a generic picture illustrating
what we have built up so far.}{\hsize}{fig:partial1}

Now $E\|2$, so that $T_E$ must overlap $T_2$.  To show that it must
occur as in Figure \ref{fig:jaws}, we assume on the contrary that
$r(E) < l(2)$, so that $L(2) \leq R(E)$.

However, $E\prec 3$, so $R(E) < L(3) < L(C)$.  We also have $E\|C$, so
that $l(C)\leq r(E)$, as shown in Figure \ref{fig:back1}.
\fig{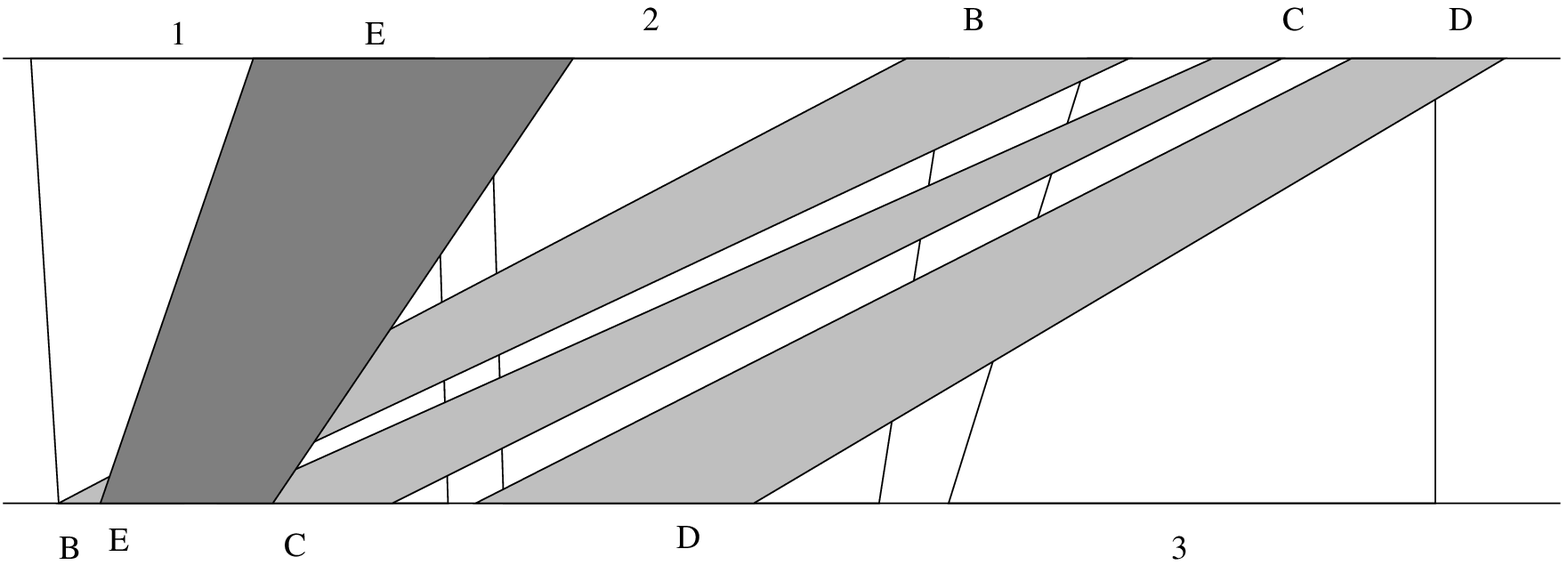}{Generic representation after inserting $E$
according to our (wrong) assumption}{\hsize}{fig:back1}

Putting the chains of Figure \ref{fig:jaworder} together with our
assumption that $L(2)\leq R(E)$, we get that $R(1) < L(2)\leq R(E)<L(F)<L(G)$.
Since $G\|1$ in the original order, this chain of inequalities,
together with the chains from the original order, imply that
$r(F)<l(G)\leq r(1)<l(D)$.  Putting this together with the fact that
$F\|D$, we find that $L(D)\leq R(F)$.

\fign{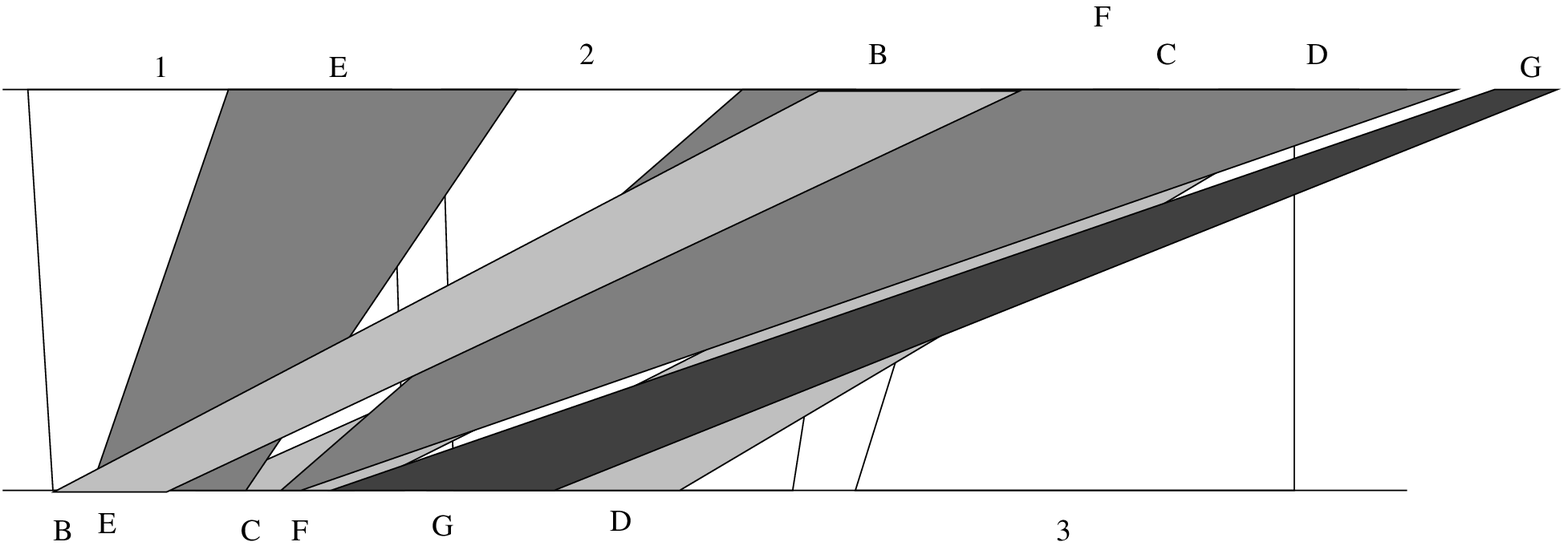}{Here, $B\prec G$, although in Figure
\ref{fig:jaworder}, $B\|G$.  $C$ and $D$ are hard to see because they
are covered by $F$ and $G$.}{\hsize}

Now everything falls apart.
Figure \ref{fig:jaworder} and $L(D)\leq R(F)$ together imply that
$R(B)<L(D)\leq R(F)<L(G)$.  Earlier we had observed that $l(C)\leq r(E)$; this
implies that $r(B)<l(C)\leq r(E)$, but $r(E) < l(G)$ (both steps from
Figure \ref{fig:jaworder}), so that $r(B) < l(G)$.  Thus, in the order
represented by these trapezoids, $B \prec G$, although in the jaw
order of Figure \ref{fig:jaworder}, $G\|B$.

Therefore we conclude that the assumption that $r(E) < l(2)$ was
wrong, and that actually
\begin{equation}\label{eq:makejaw}r(E) \geq l(2)\end{equation}
giving us
$$r(B)\arlt l(C)\overset{\eqref{eq:Cleft}}{\arlte} r(1)\arlt l(2)
\overset{\eqref{eq:makejaw}}{\arlte} r(E)\arlt
l(D)\overset{\eqref{eq:assumption}}{\arlte} r(2)\arlt l(3)\leq r(F) \arlt
l(G)$$

However, $E\prec D$, so that forces the lower left endpoint of $D$
farther to the right than our picture in Figure \ref{fig:partial1}.

There are two inequalities remaining to be shown.  Since we have just
established that $l(2)\leq r(E)$, both $F$ and $G$ must be to the right of
1 on the lower baseline, and hence, to the left of or overlapping on
the upper baseline.  Therefore, $L(F)<R(1)$, and so because $E \prec
F$ and $1 \prec 2$,
\begin{equation}\label{eq:jawout}R(E)<L(2)\end{equation}

$$R(E)\overset{\eqref{eq:jawout}}{\arlt} L(2)\arlte R(2)
\overset{\eqref{eq:incomp1}}{\arlt} L(D)$$

One to go.  We have just observed that $G$ overlaps 1 on the upper
baseline.  Thus, since $F \prec G$, $L(F)\leq R(F)<L(G)\leq R(1)$.  $F$ must 
overlap 3, so the lower interval for $F$ must appear to the right of
or overlap the lower interval for 3.  i.e.,
\begin{equation}\label{eq:lastone}l(3)\leq r(F)\end{equation}
$$ r(B)\arlt l(C)\overset{\eqref{eq:Cleft}}{\arlte} r(1)\arlt l(2)
\overset{\eqref{eq:makejaw}}{\arlte} r(E)\arlt l(D)
\overset{\eqref{eq:assumption}}{\arlte} r(2)\arlt l(3)
\overset{\eqref{eq:lastone}}{\arlte} r(F)\arlt l(G)$$
\end{proof}

The lemma asserts a long chain of inequalities on endpoints.  A
representation that respects all of these inequalities is shown in
Figure \ref{fig:finaljaw}.  Most of these inequalities are only
important for technical reasons.  The key idea that we will be
exploiting is the chains of inequalities pictured in Figure
\ref{fig:jaws}.  The relative positions of the endpoints in this piece
of the picture are forced by the inequalities, though there is room
for movement in the other areas of the picture, as one can see by
moving endpoints of other intervals in Figure \ref{fig:finaljaw}.

\fig{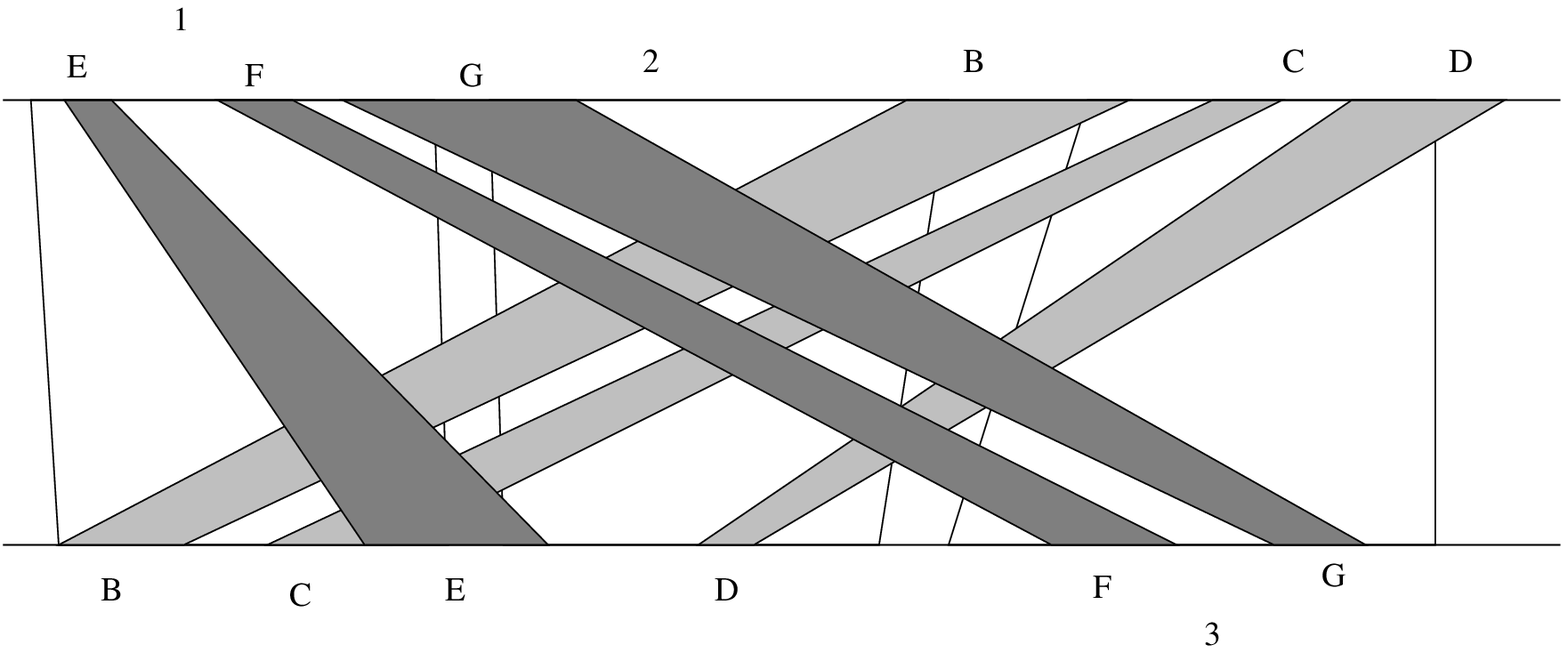}{Trapezoid representation of the jaw order}
{\hsize}{fig:finaljaw}

\section*{Improper Trapezoid Orders}
Now, with this lemma in hand, we proceed to produce an example of an
improper trapezoid order.

\begin{thm}
\label{thm:improper}
The order shown in Figure \ref{fig:improper} is an improper trapezoid
order.  \fig{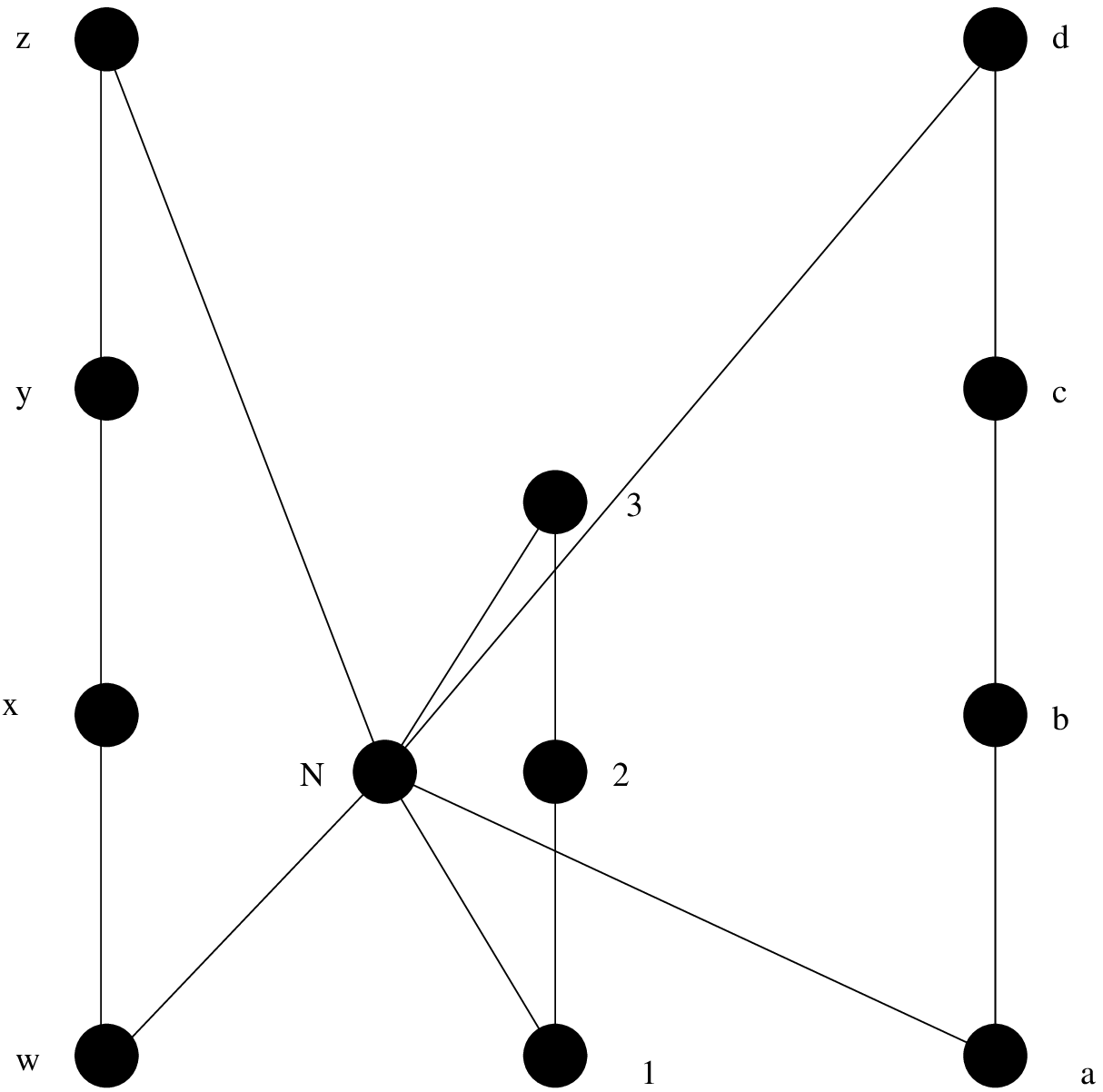}{An improper trapezoid
order}{3in}{fig:improper}
\end{thm}

\fig{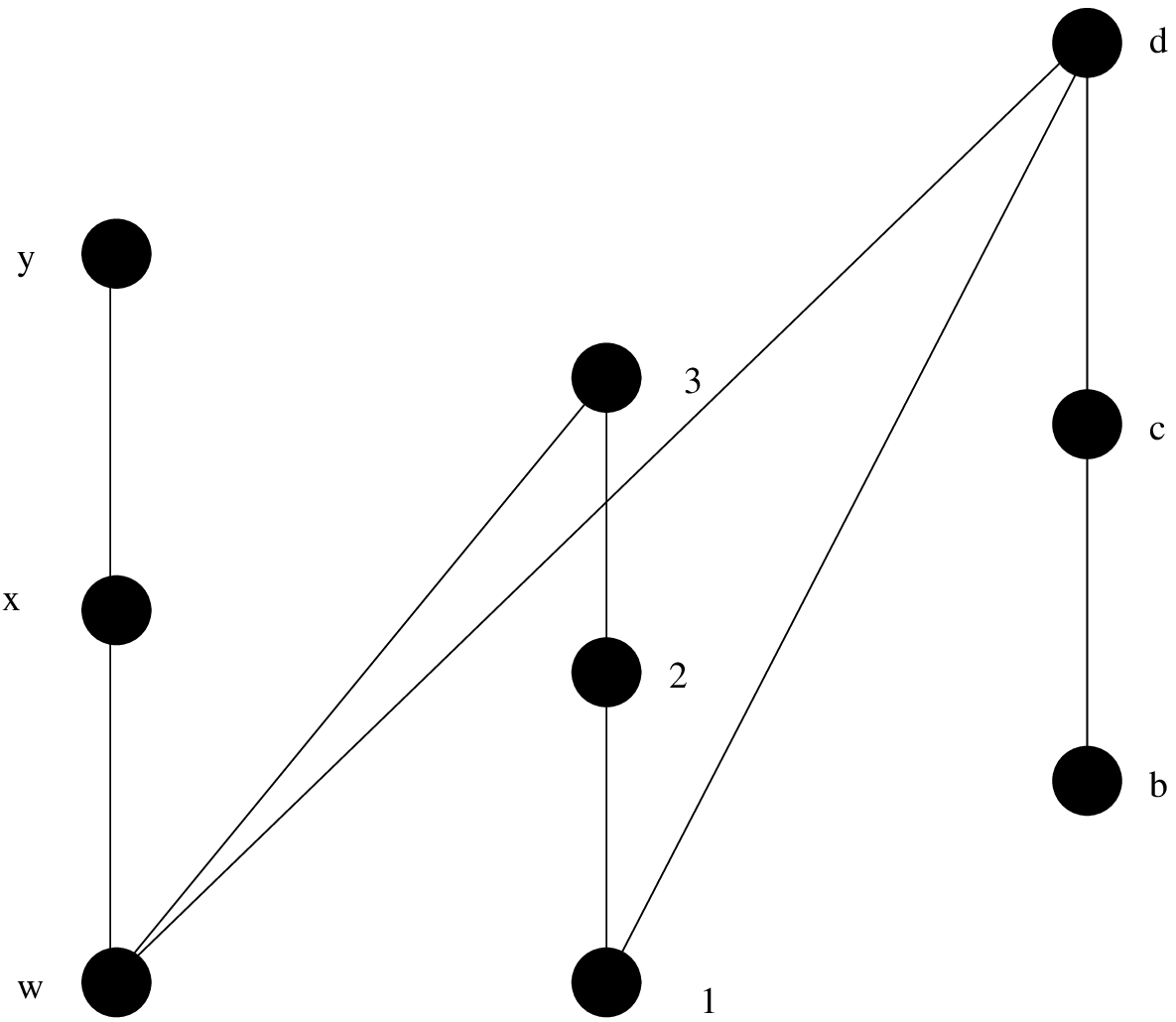}{One jaw order inside Figure \ref{fig:improper}}{3in}
{fig:jawrestriction}
\fig{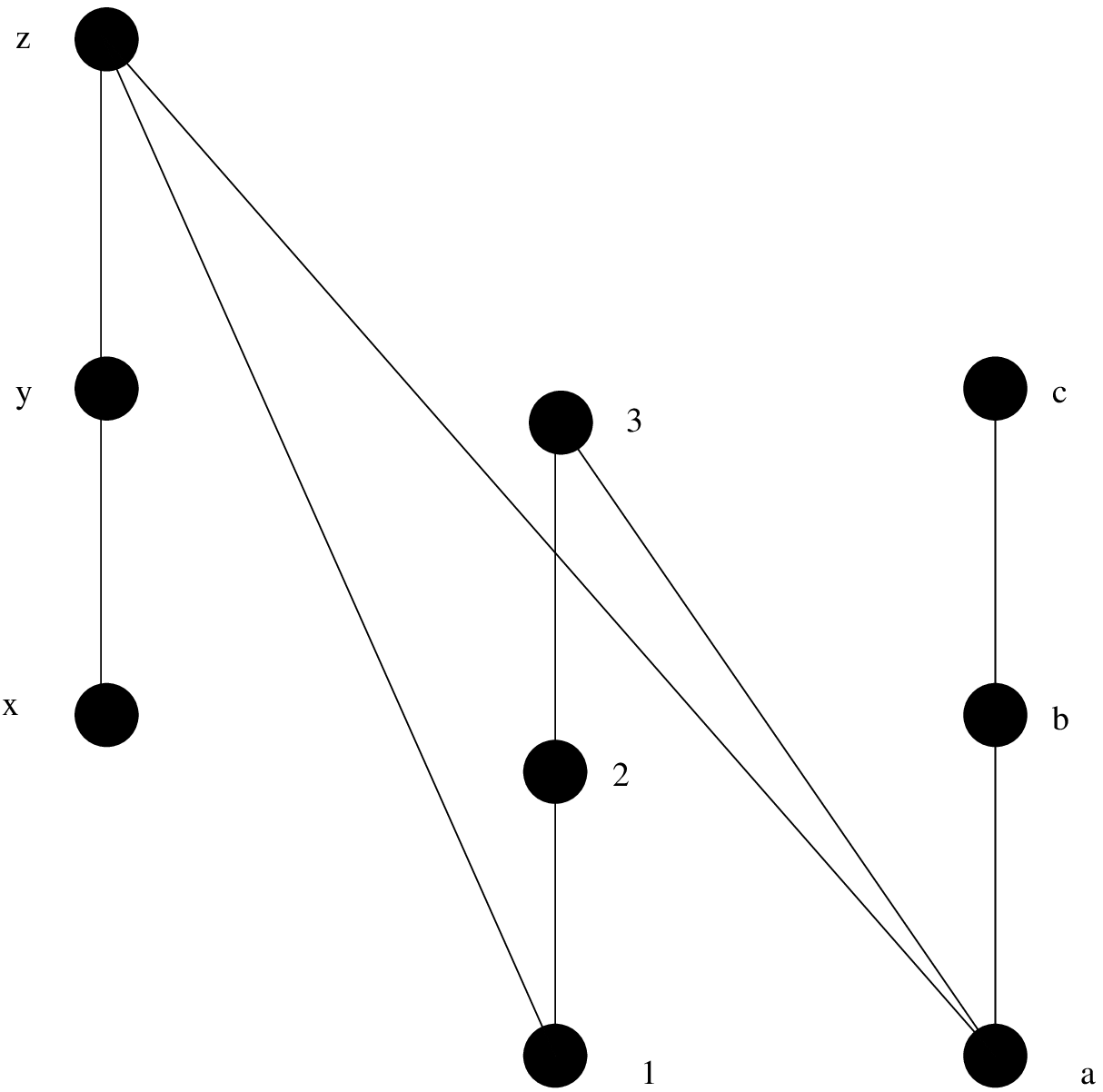}{Another jaw order inside Figure \ref{fig:improper}}{3in}
{fig:dualjawrestriction}
\begin{proof}
There are two restrictions of this order which are isomorphic to the
jaw order.  By the jaw lemma, each of these restrictions must produce
a jaw formation like Figure \ref{fig:jaws} in any trapezoid representation.

Let us temporarily suppose that both jaw formations are on the lower
baseline.  From applying the jaw lemma to Figure
\ref{fig:jawrestriction}, we get
$$r(b)<l(c)\leq r(1)<l(2)\leq r(w)<l(d)\leq r(2)<l(3)\leq r(x)<l(y)$$
From applying the jaw lemma to Figure \ref{fig:dualjawrestriction}, we
get
$$r(x)<l(y)\leq r(1)<l(2)\leq r(a)<l(z)\leq r(2)<l(3)\leq r(b)<l(c)$$
Combining these two, we find that in the first, $l(3)\leq r(x)$, but in
the second, $r(x)<l(3)$.  We conclude that our temporary supposition
was wrong, and that the two jaw formations must be on opposite
baselines.

\fig{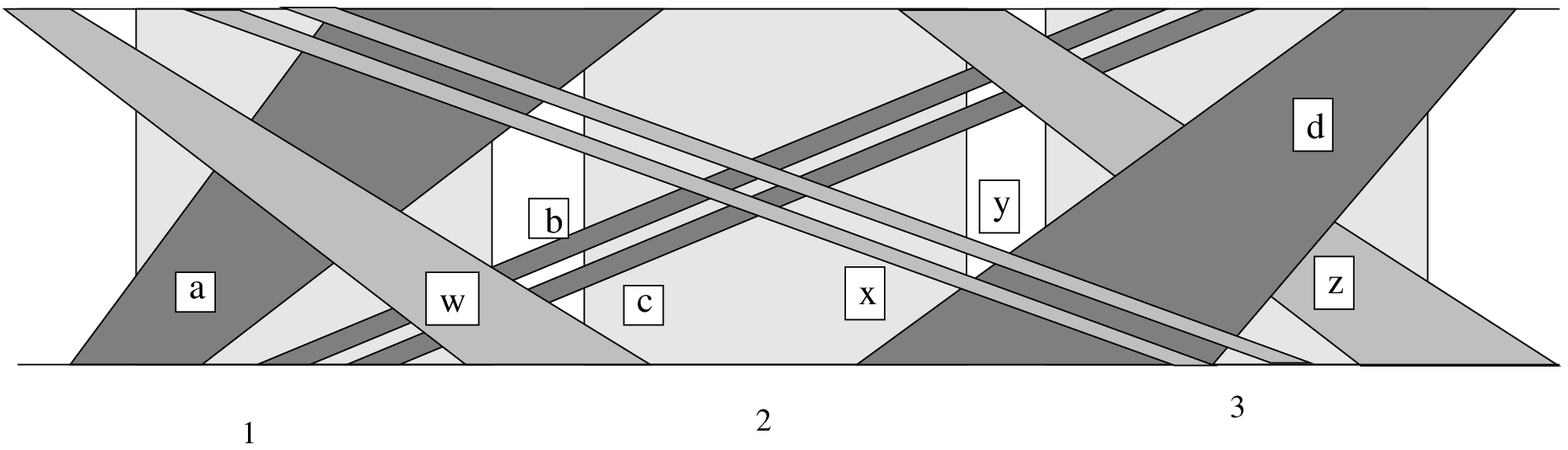}{The framework built by two applications of the Jaw Lemma}
{\hsize}{fig:frame}
Thus, we obtain the representation shown in Figure \ref{fig:frame},
with $l(2) \leq r(w) < l(d) \leq r(2)$ and $L(2)\leq R(a)<L(z)\leq
R(2)$.  Now we look for the place to put in $N$.  Since $N$ is clamped
by the jaws on both sides -- i.e. since $a \prec N \prec z$ and $w
\prec N \prec d$ -- the intervals for $N$ must be completely contained
in the intervals for $2$, and hence any trapezoid representation of
this order must be improper, as the representation given in Figure
\ref{fig:improptrap}.  (Note that the representation in Figure
\ref{fig:improptrap} could have been drawn completely using only
parallelograms, and hence this is also an example of an improper
parallelogram order.)

\fig{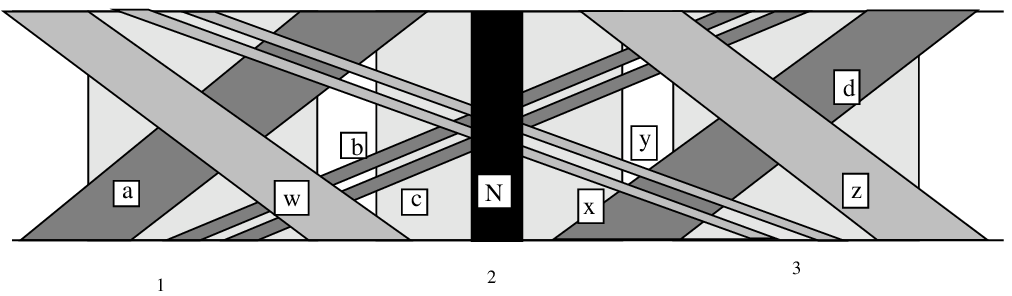}{Improper parallelogram representation of
the order in Figure \ref{fig:improper}}{\hsize}{fig:improptrap}

\end{proof}

We now turn our attention to showing the existence of a trapezoid
graph with no proper representation.  Gallai \cite{gallai67:_trans} 
shows that two orders have the same comparability graph if and
only if there exists a finite sequence of orders such that each is
obtained from the previous by reversing an autonomous set of the
order.  An autonomous set is a subset of the elements in the order
such that each element not in the subset has the same relation to all
elements of the subset.  i.e. $A$ is an autonomous set of
$P=(X,\prec)$ if $A \subseteq X$ and $\forall t \not \in A$, $t \prec
x$ for any $x \in A$ implies that $t \prec y$ for all $y \in A$, $t
\succ x$ for any $x \in A$ implies that $t \succ y$ for all $y \in A$,
and $t \| x$ for any $x \in A$ implies that $t \| y$ for all $y \in
A$.

It therefore suffices to prove that a property is invariant with
respect to reversing autonomous sets of an order to prove that it is
true of the corresponding comparability graph (or incomparability
graph, as the case may be.)

A simple observation from the definition of an autonomous set is that
if $a,b,c \in X$, $A$ is an autonomous set, and $a,b \in A$, then if
$a \prec c$ (or $a \succ c$) and $b\|c$, then it must be the case that
$c \in A$ as well.  Otherwise, $A$ would not be an autonomous set.

\begin{prop}\label{prop:autonomous}
The only non-trivial autonomous sets in the order of Figure
\ref{fig:improper} are $\{b,c\}$ and $\{x,y\}$.
\end{prop}

\begin{proof}
The given sets are clearly autonomous sets, so we only need to prove
that no other proper subsets of at least two elements can be
autonomous.  (Singletons are automatically autonomous sets, but not
very interesting ones.)

Suppose that $A$ is a non-trivial autonomous set such that $a \in A$,
and also suppose that $b \not \in A$.  Then, since $b \succ a$,
$\forall t \in A, b \succ t$, i.e. $A \subseteq Pred(b)$.  However,
$Pred(b) = \{a\}$, contradicting non-triviality of $A$.  Thus, $b \in
A$.  Now, $\forall t \in \{3,z,N\}, t \succ a$ and $t \| b$.
Therefore, it must be the case that $t \in A$ as well.  Since $z \in
A$, by symmetry we must also have that $\{y,1,a,N\} \subseteq A$.
Once again, $w \prec N$ and $w \| b$ imply that $w \in A$, and so
again by symmetry, $\{x,c,d\} \subseteq A$.  Finally, $2 \prec 3$ and
$2 \| N$ imply that $2 \in A$, so that $A$ is the whole order.

By symmetry, if $A$ is any nontrivial autonomous set containing any of
$\{a,d,w,z\}$, then $A$ must be the whole order.

By analogous arguments, if $A$ is any nontrivial autonomous set
containing $1$ or $3$, then it must be the whole order.

Thus, the only possible nontrivial autonomous sets are those
containing only $\{b,c,x,y,2,N\}$.  If $A$ is some autonomous set
containing 2, then since $1 \prec 2$, $1 \prec t$ for every $t \in A$,
or equivalently, $A \subseteq Succ(1)$.  The only possible element for
$A$ other than 2 in this case is $N$; however, $a \prec N$ and $a \|
2$, so $\{2,N\}$ is not an autonomous set.  Similarly, if $N \in A$,
$1 \prec N$ but $1 \| \{b,c,x,y\}$, so that would imply that $1 \in
A$, a case we already discovered would mean $A$ was the whole order.

Thus, any nontrivial autonomous set must contain only elements from
$\{b,c,x,y\}$.  Since $\{b,c\}$ and $\{x,y\}$ are known to be
autonomous sets, we only need to consider what happens if, say $b \in
A$ and $x$ or $y$ is in $A$.  In either of these cases, though, $z \|
b$ and $z \succ x$ and $z \succ y$, so that $z \in A$ and again, $A$
must be the whole order.  By symmetry, the analogous case involving
$c$ instead of $b$ is also covered.
\end{proof}

\begin{prop}\label{prop:compinvariant}
The trapezoid graph represented by Figure \ref{fig:improptrap} is an
improper trapezoid graph.
\end{prop}

\begin{proof}
Since the order in Figure \ref{fig:improper} has only two non-trivial
autonomous sets, and we know that any trapezoid representation of that
order must be an improper trapezoid representation, we only need to
check that reversing the given autonomous sets preserves the improper
trapezoid representation.  However, simply reversing the trapezoids
for $b$ and $c$ or for $x$ and $y$ doesn't change the picture at all,
only the labels, and so the result must also be an improper trapezoid
order; hence, the associated cocomparability graph must be an improper
trapezoid graph.
\end{proof}

\section*{Unit Trapezoid Orders}
A standard observation is that a unit interval order must be a proper
interval order, because there is no way for one interval to be
properly contained in another, yet have the same length.  Similarly,
every unit trapezoid order must be a proper trapezoid order, because
there is no way for one trapezoid to be properly contained in another
if they both have the same area.  In the case of interval orders,
every proper interval order is also a unit interval order; i.e. the
two classes of orders are the same.  We now proceed to show that this
is not the case for trapezoid orders.

\fig{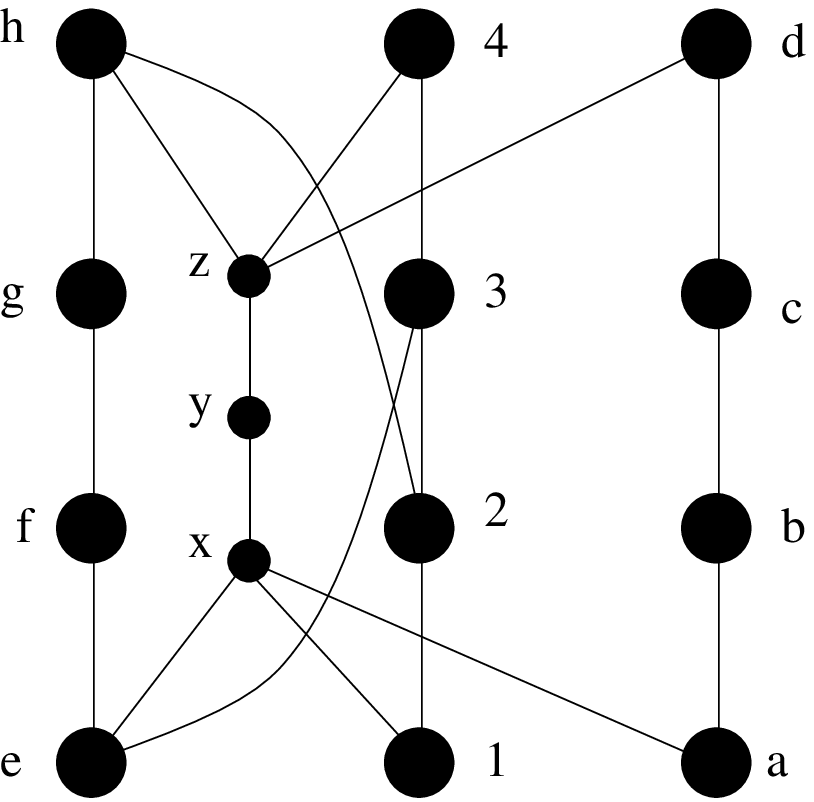}{An example of a proper trapezoid order with no unit 
trapezoid representation.}{3in}{fig:pnuorder}

\begin{thm}\label{thm:pnu}
The order given in Figure \ref{fig:pnuorder} is a proper
trapezoid order with no unit trapezoid representation.
\end{thm}

\begin{proof}
To show that the order in Figure \ref{fig:pnuorder} is a proper
trapezoid order, it suffices to exhibit a proper trapezoid
representation.  This is provided in Figure \ref{fig:pnutrap}.  Notice
that this, too, can be drawn using parallelograms, and hence the order
in Figure \ref{fig:pnuorder} is a proper parallelogram order.
\fig{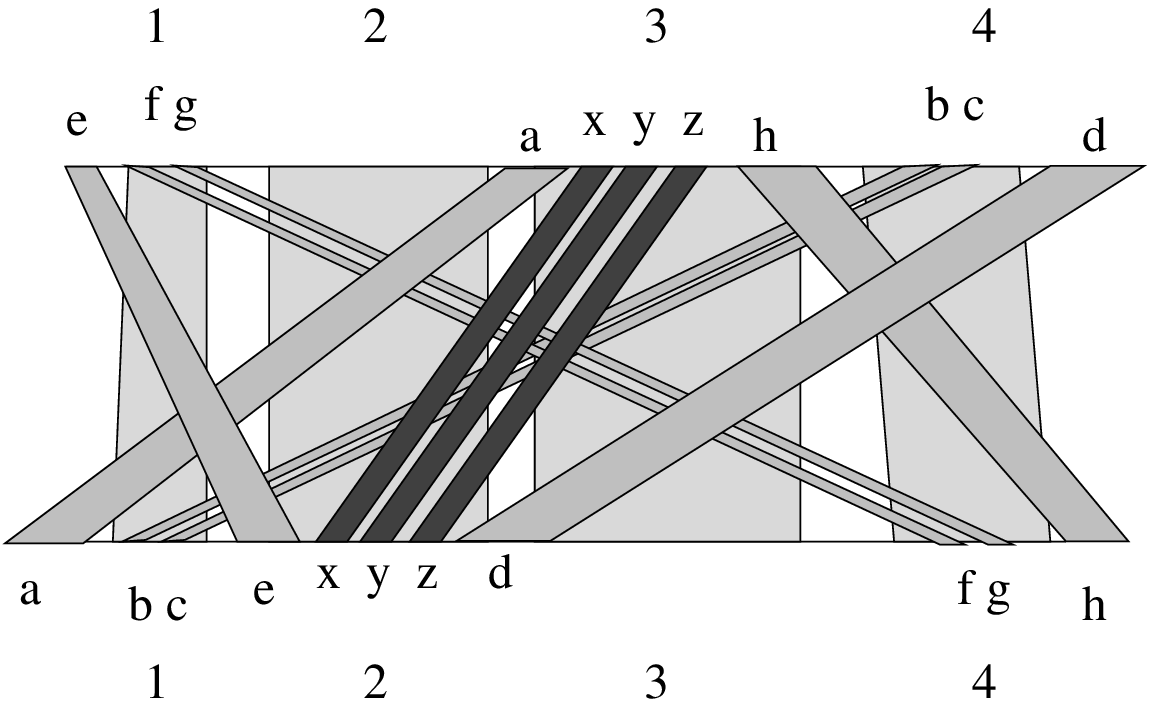}{A proper trapezoid representation of Figure 
\ref{fig:pnuorder}}{\hsize}{fig:pnutrap}

To show that it has no unit trapezoid representation, we will again
make use of the jaw lemma.


\begin{figure}[ht]
{\noindent
\hfill
\begin{minipage}[b]{.4\linewidth}
\centering\epsfig{figure=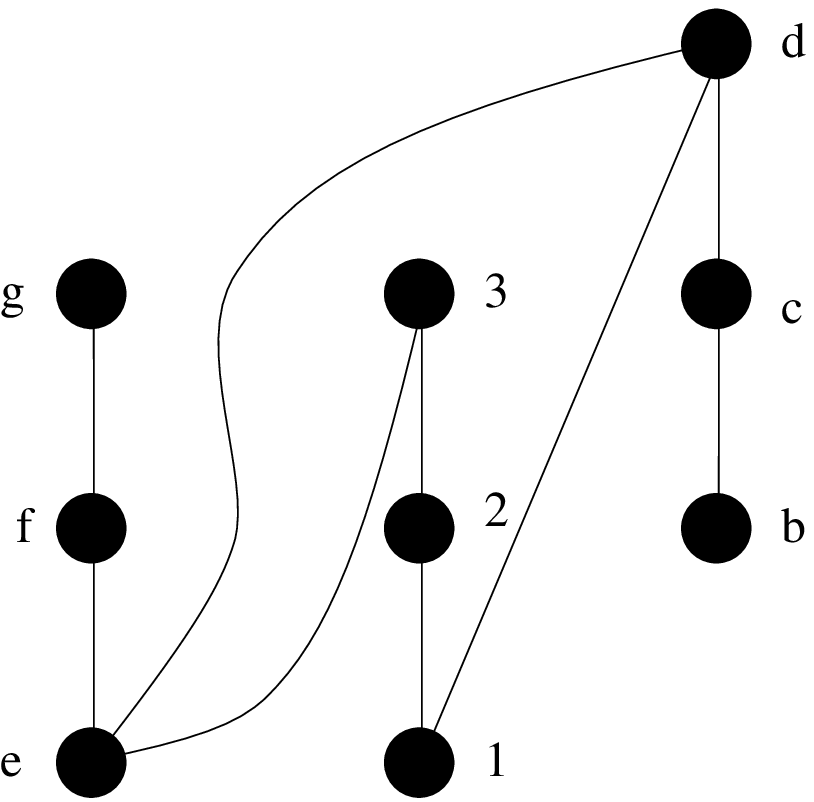,width=\linewidth}
\end{minipage}\hfill
\begin{minipage}[b]{.4\linewidth}
\centering\epsfig{figure=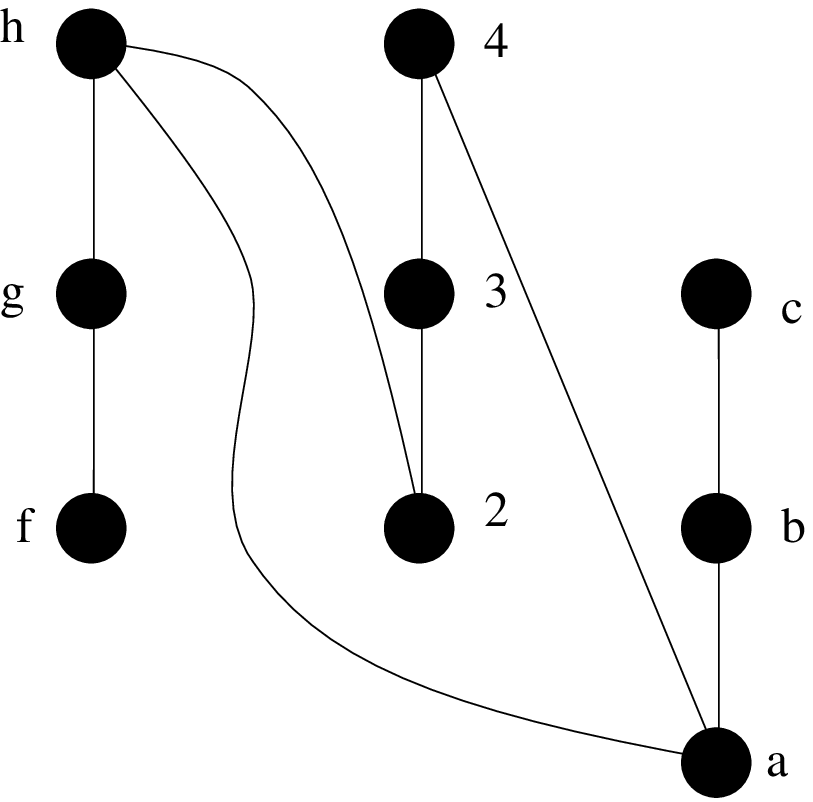,width=\linewidth}
\end{minipage}
\hfill
\caption{Two restrictions of Figure \ref{fig:pnuorder} isomorphic to Figure 
\ref{fig:jaworder}}
\label{fig:pnurest}
}
\end{figure}

The first order in Figure \ref{fig:pnurest}, by the jaw lemma, must
appear in a jaw formation, with the trapezoids for $d$ and $e$ as the
``teeth''.  Without loss of generality, assume that the formation
appears on the lower baseline.

The second order in Figure \ref{fig:pnurest} is also the jaw order,
and hence we can also apply Lemma \ref{lem:jaw} to it to get that it
too must appear in a jaw formation.  As with the preceding example,
this formation must appear on the upper baseline.  The proof is the
same as the proof in the preceding example, with only small
modifications to accomodate the slightly different restrictions.

The result of this is that the restriction of any trapezoid
representation to the elements in Figure \ref{fig:pnurest} must be
similar to the one given in Figure \ref{fig:pnuframetrap}.  In
particular, the trapezoids for 2 and 3 and the overlaps of other
trapezoids with those for 2 and 3 must appear as illustrated.

\fig{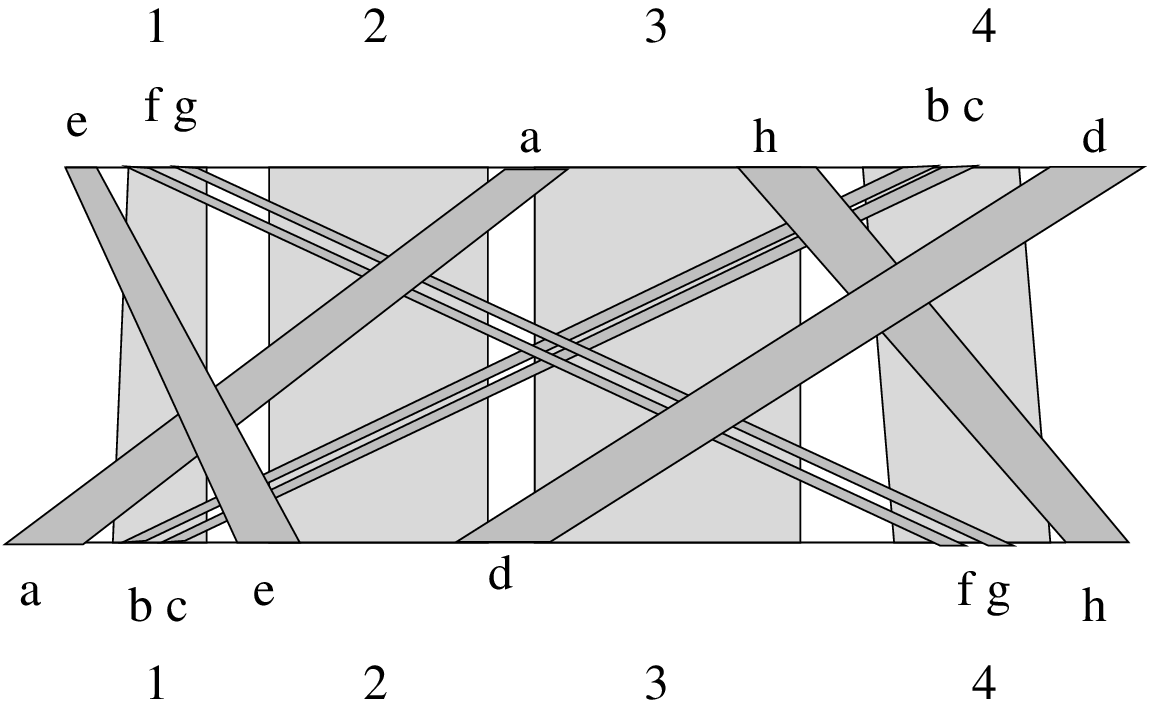}{This framework is forced by the jaw lemma; 
the representation in Figure \ref{fig:pnutrap} follows easily.}{\hsize}
{fig:pnuframetrap}

Now there is only one place to insert the chain $x\prec y\prec z$,
namely between the jaws on each side, resulting in a representation
similar to that in Figure \ref{fig:pnutrap}.

We now conclude that no trapezoid representation of this order can be
a constant area representation.  The area of any trapezoid is given by
the formula $\frac 1 2 (b_1 + b_2) h$.  For all of the trapezoids in
any representation, the height and the $\frac 1 2$ are constant, so
that for the trapezoids to have constant area, the sum of the bases
must be constant.  By a suitable choice of scaling, we may assume that
the sum of the bases is 2.  This implies that $\left [ t(x)+t(y)+t(z)
\right ] + \left [ b(x)+b(y)+b(z) \right ] = 6$.  Hence, at least one
of the two summands must be at least 3. However, the intervals for
$x$, $y$ and $z$ are all contained inside the lower interval for 2 and
the upper interval for 3.  Thus, either the lower interval for 2 or
the upper interval for 3 would have length strictly greater than 3,
which would mean that the sum of the lengths for the intervals of at
least one of the two is not 2; thus, the representation cannot be a
constant area representation.

This example is not minimal; for example, if we remove $y$ from this
order, the resulting restriction is still a proper trapezoid order
with no unit trapezoid representation; however, the proof is slightly
less obvious.  There may also be smaller examples.
\end{proof}

\begin{prop}The only non-trivial autonomous sets in Figure
\ref{fig:pnuorder} are $\{b,c\}$, $\{f,g\}$, $\{x,y\}$, $\{y,z\}$, and
$\{x,y,z\}$.
\end{prop}

\begin{proof}The arguments here are entirely analogous to those in
Proposition \ref{prop:autonomous}, though significantly longer (and
more boring) due to the loss of some symmetry, and so will not be
repeated here.
\end{proof}

\begin{prop}The trapezoid graph represented in Figure
\ref{fig:pnutrap} is a proper trapezoid graph, but not a unit
trapezoid graph.
\end{prop}

\begin{proof} Almost identical to the proof of proposition
\ref{prop:compinvariant}.
\end{proof}

The arguments used to show that our examples give rise to examples of
improper trapezoid graphs and proper but not unit trapezoid graphs
were ad hoc.  A more satisfying approach would be to show that the
properties of being a proper trapezoid order and being a unit
trapezoid order are comparability invariants in the sense of
\cite{habib91:_inter} and \cite{felsner:_semi}.  However, the
techniques of \cite{habib91:_inter} and \cite{felsner:_semi} do not
seem to apply directly to proper and unit trapezoid orders, leading us
to the ad hoc arguments we used.

\end{document}